\documentclass[a4paper,12pt]{amsart}
\usepackage{amsmath,
			amssymb,
			amsthm,
			amsfonts,
			amscd,% Creo que no lo vamos a usar mucho, es provisional
			mathrsfs, % letras caligráficas bonitas de foliaciones
			enumerate, % para poner [i.] y todo eso
			mathtools} % Para \mathclap http://tex.stackexchange.com/questions/12344/separate-long-math-text-under-sum-symbol-into-different-lines/33358#33358
\usepackage[utf8]{inputenc}
\usepackage{psfrag}% I added them
\usepackage[isbn=false,doi=false,eprint=false,style=alphabetic]{biblatex}   % hay que usar biber en vez de bibtex
\usepackage{graphicx}
\usepackage[all]{xy}
\usepackage[normalem]{ulem}
\usepackage{url} % para que bibtex ponga las url.
\usepackage{color}
\usepackage[top=1.5in,
			bottom=1.5in,
			left=1in,
			right=1in]{geometry}
\usepackage{braket} % para usar \set y \Set
\usepackage{fontawesome5}

\usepackage{etoolbox}
\makeatletter
\patchcmd{\@maketitle}
  {\ifx\@empty\@dedicatory}
  {\ifx\@empty\@date \else {\vskip3ex \centering\footnotesize\@date\par\vskip1ex}\fi
   \ifx\@empty\@dedicatory}
  {}{}
\patchcmd{\@adminfootnotes}
  {\ifx\@empty\@date\else \@footnotetext{\@setdate}\fi}
  {}{}{}
\makeatother
\newdimen\theight
\def\TeXref#1{%
             \leavevmode\vadjust{\setbox0=\hbox{{\tt
                     \quad\quad  {\small \textrm #1}}}%
             \theight=\ht0
             \advance\theight by \lineskip
             \kern -\theight \vbox to
             \theight{\rightline{\rlap{\box0}}%
             \vss}%
             }}%

\newtheorem{theorem}{Theorem}[section]

\newtheorem{proposition}[theorem]{Proposition}
\newtheorem{corollary}[theorem]{Corollary}
\theoremstyle{definition}
\newtheorem{definition}[theorem]{Definition}
\newtheorem{example}[theorem]{Example}

\theoremstyle{remark}
\newtheorem{remark}[theorem]{Remark}

\numberwithin{equation}{section}

\newcommand\chii{\raise2pt\hbox{$\chi$}}
\newcommand\phii{{\raise2pt\hbox{$\varphi$}}}

\newcommand\ese{\mathbb{S}}

\newcommand{\R}{\mathbb{R}}

\newcommand\F{\mathscr{F}}

\newcommand\LL{\mathscr{L}}

\newcommand{\THL}[1]{(THL)_{#1}}
\newcommand{\HL}[1]{(HL)_{#1}}
\newcommand{\HB}[1]{H_B^{#1}}
\newcommand{\HM}[1]{H_M^{#1}}
\newcommand{\PHB}[1]{P\!H_B^{#1}}
\title{Hard Lefschetz Property for $\ese^3$-actions}

\author[J.I.~Royo Prieto]{Jos\'{e} Ignacio Royo Prieto}
\address{Matematika Saila\\ Zientzia eta Teknologia Fakultatea\\ University of the Basque Country UPV/EHU\\ Barrio Sarriena s/n\\ 48940 Leioa\\Spain.}
\email{joseignacio.royo@ehu.eus}
\thanks{Partially supported by Ministerio de Ciencia, Spain, grant PID2019-105621GB-I00. The authors acknowledge that the research cooperation was funded by the program Excellence Initiative Research University at the Jagiellonian University in Krakow within the framework of the research
	group Reeb-Reinhart 2022.}

\author[M.~Saralegi Aranguren]{Martintxo Saralegi-Aranguren}
\address{F\`{e}d\`{e}ration CNRS\\  Nord-Pas-de-Calais FR 2956\\ UPRES-EA 2462 LML\\
Facult\'e Jean Perrin\\ Universit\'{e} d'Artois\\   Rue Jean Souvraz SP 18\\ 62 307 Lens Cedex, France.}
\email{martin.saraleguiaranguren@univ-artois.fr}

\author[R.~Wolak]{Robert Wolak}
\address{Instytut Matematyki\\ Uniwersytet Jagiellonski\\
	ul. prof. Stanis{\l}awa {\L}ojasiewicza 6
	30-348 Krak\'ow,  Poland.}
\email{robert.wolak@im.uj.edu.pl}

\keywords{Lefschetz Hard Property, 3-Sasakian Manifolds, Gysin sequence, $\ese^3$-actions}

\subjclass[2010]{53C12, 53D10, 53C25}

\date{\today}

\begin{document}

\begin{abstract}
The Hard Lefschetz Property (HLP) has recently been formulated in the context of isometric flows without singularities on manifolds. In this category, two versions of the HLP (transverse and not) have been proven to be equivalent, thus generalizing what happens in the important cases of both K-contact and Sasakian manifolds. In this work we define both versions of the HLP for almost-free $\ese^3$-actions, and prove that they agree for actions satisfying a cohomological condition, which includes the important category of 3-Sasakian manifolds, where those two versions of the HLP are shown to be held. We also provide a family of examples of free actions of the 3-sphere which are not 3-Sasakian manifolds, but satisfy the HLP.
\end{abstract}

\maketitle

\section*{Introduction}

The Hard Lefschetz Property (HLP) is a relevant cohomological duality property which has traditionally been studied in several categories of the symplectic world  (we refer to the introduction of \cite{leches} for a brief account of them). That duality is to be satisfied (or not) between cohomology groups of complementary dimensions by means of repeated product by an special 2-class, whose representative is the symplectic form (in the context of symplectic manifolds), a contact form (contact or Sasakian manifolds) or the Euler form (isometric flows). 

For Sasakian manifolds (which are a special case of both contact manifolds and isometric flows) this duality is satisfied by the basic cohomology of the manifold, so the HLP is a {\em transverse} property. In \cite{mino} the authors give a characterization of the HLP satisfied by the Sasakian manifolds in terms of isomorphisms of the cohomology groups of the manifold itself (not basic), which leads them to define the Lefschetz property for any contact manifold in the same global terms. Both definitions of the HLP, transversal and not, were proven to be equivalent by \cite{linyi} in the category of K-contact flows (contact manifolds which are, at the same time, isometric flows). 

We have recently shown in \cite{leches} that the HLP is not exclusive of the symplectic framework by formulating it for any isometric flow (or, equivalently, any Killing vector field without singularities) and proving that both versions, transversal and not, of the HLP are equivalent. These definitions of the HLP agree with those already given for K-contact flows. This highlights the cohomological nature of the HLP. The Gysin sequence for isometric flows is an essential tool for the construction, for its connecting morphism is the multiplication by the Euler class, which is the basic 2-class involved in the definition of the HLP.

It is natural to extend the definition of the HLP to the framework of almost-free actions of the 3-sphere $\ese^3$, a context where we also have a Gysin sequence (it is a classical tool when we have a free $\ese^3$-action, but it also works for general smooth actions of $\ese^3$, see \cite{debrecen}). In this paper we define the HLP for this kind of actions both for the basic cohomology (the isomorphisms being realized by repeated multiplication of the corresponding Euler class, which is a basic 4-class associated to the action) and for the cohomology of the manifold itself. The main result of Section 2 shows that both versions of the HLP agree when a cohomological requirement is satisfied: the multiplication by the Euler class must induce an injective map in the basic cohomology. We call it the {\em injectivity condition}.

3-Sasakian manifolds have been intensively studied since the 60's of the last century. Due to interesting applications in physics, the recent years have witnessed a renewed interest in these structures (see \cite{BG99} and \cite{galicki} for a detailed review). 3-Sasakian manifolds naturally carry an almost-free action of $\ese^3$. In Section 3, we show that any Sasakian manifold satisfies the injectivity condition. We also show that any 3-Sasakian manifold satisfies the HLP, which was already known for the transversal version of the HLP when the sectional curvature is positive (see \cite{mohseni}).

%In this work we define Lefschetz property HLP for $\ese^3$-actions in local and global fashion, and show its cohomological nature by showing the equivalence. We also prove that 3-Sasakian manifolds satisfy both properties. 

We finish the paper by providing a family of $\ese^3$-actions on manifolds of dimension $4n+3$ for any $n\ge 2$ which satisfy the HLP, but fail to admit compatible 3-Sasakian structures, showing that the HLP is not exclusive of 3-Sasakian manifolds.

\section{Almost free actions of the 3-sphere}

%\begin{itemize}
%	\item Gysin sequence
%	
%	\item Induced $\ese^1$-action and Gysin sequence with $B$ and $Z$.
%	
%	\item Injectivity $H^*_B\to H^*_Z$.
%
%	\item $[e_1^2]_Z=[\Omega^2]_B$
%
%\end{itemize}
%
%For this brief version we shall not go into detail with this.

We refer the reader to \cite{debrecen} and \cite{leches} for the basics regarding actions of $\ese^3$ on manifolds and isometric flows, respectively. Let $M$ be a closed manifold of dimension $4n+3$, and consider a smooth almost free $\ese^3$-action $\Phi$ on $M$. Recall that the orbit space $B=M/\ese^3$ is an orbifold of dimension $4n$. The $\ese^3$-action induces a foliation of dimension $3$ on $M$, which we denote by $\F_{\ese^3}$. For short, we denote by $H_B^*=H^*(M/\F_{\ese^3})$ its basic cohomology, which is the cohomology of the differential complex of forms of $M$ which are both invariant and horizontal for the action:
$$
\Omega_B^*=\set{\omega\in\Omega^*_M | i_X\omega = L_X\omega =0 
	\text{ for any } X\in T\F_{\ese^3}}.
$$

{\bf Notation regarding the action:} Denote by $X_i,\ i=1,2,3$ the three fundamental vector fields of the $\ese^3$-action. Notice that each $X_i$ is the fundamental vector field of an almost free $\ese^1$-action on $M$. Taking an $\ese^3$-invariant metric $\mu$, the dual characteristic forms are $\chi_i=i_{X_i}\mu$. We have
$$
\begin{aligned}
d\chi_1&=e_1 - \chi_2\wedge \chi_3\\
d\chi_2&=e_2 + \chi_1\wedge \chi_3\\
d\chi_3&=e_3 - \chi_1\wedge \chi_2	
\end{aligned}
$$
where $e_i\in\Omega^2(M)$ are ${\ese^3}$-horizontal forms. Write $\F_i$ the foliation induced by $X_i$ (which are, indeed, isometric flows). Notice that $e_i$ is $\F_i$-basic and is, indeed, the Euler form for $\F_i$, but $e_i$ is not $\ese^3$-basic. Nevertheless, it is straightforward to check that
\begin{equation}\label{eq:omega}
d(\stackrel{\Theta}{\overbrace{e_1\wedge\chi_1 + e_2\wedge\chi_2 + e_3\wedge\chi_3}})=
\stackrel{\Omega}{\overbrace{e_1^2 + e_2^2 + e_3^2}} + 
d(\stackrel{\Upsilon}{\overbrace{\chi_1\wedge\chi_2\wedge\chi_3}}).
\end{equation}
 The 4-form $\Omega=e_1^2 + e_2^2 + e_3^2$ is basic and we call it the {\em Euler form} of the $\ese^3$-action. We call the class $[\Omega]\in H^4_B$ the {\em Euler class} of $\Phi$ (or, for short, the Euler class, if no confussion appears).

\begin{remark}\label{rem:su2}
For the almost-free $\ese^3$-action arising from a 3-Sasakian structure, in \cite[formula~(13.5.2)]{galicki} the identity \eqref{eq:omega} appears as $d\Theta=\Omega+2d\Upsilon$, which is due to the fact that in that book the fundamental fields, denoted as $\set{\xi_1,\xi_2,\xi_3}$, satisfy $[\xi_a,\xi_b]=2\epsilon_{abc}\xi_c$, while we have chosen the convention $[X_a,X_b]=\epsilon_{abc}X_c$ for the Lie algebra $su(2)$ instead, accordingly to \cite{debrecen}. The metric, the fundamental forms and  the Euler forms are denoted in \cite{galicki} by $g$, $\eta_i$ and $\bar{\Phi}_i$, respectively.  By considering $\xi_i=2X_i$, and then, 
$
g=\frac{1}{4}\mu, \eta_i=\frac{1}{2}\chi_i$ and $\bar{\Phi}_i=\frac{1}{2}e_i,
$
for $i=1,2,3$, one immediately retrieves  \eqref{eq:omega} from \cite[formula~(13.5.2)]{galicki} and viceversa.

\end{remark}

If the action $\Phi$ is free, we have a principal $\ese^3$-bundle and the associated classical Gysin sequence. In \cite{debrecen} it is shown that the same exact sequence works for an almost-free action by using basic cohomology and thus avoiding to work with orbifolds, just working with forms defined in $M$. The Gysin sequence we have is:
\begin{equation}\label{eq:gysin}
\dots\to H^i_B\stackrel{\iota}{\longrightarrow} 
H^k_M \stackrel{\rho\circ\varphi}{\longrightarrow}
H^{k-3}_B\stackrel{\wedge[\Omega]}{\longrightarrow}
H^{k+1}_B  \stackrel{\iota}{\longrightarrow} 
H_M^{k+1}\to\dots
\end{equation}
where $\iota$ is induced by the inclusion of forms, $\varphi\colon H_M^k\to (H_M^k)^{\ese^3}$ is the isomorphism between the de Rham cohomology of $M$ and the  $\ese^3$-invariant cohomology,  $\rho$ is the multiple contraction $\rho([\alpha])=[i_{X_1X_2X_3}\alpha]$ and the connecting morphism is the multiplication by the Euler class.

\section{Hard Lefschetz Property}

In this section we shall define the HLP for almost-free actions of $\ese^3$. We shall mainly adapt the definitions, results and proofs of \cite{leches} from flows to our new context.

Denote by $L\colon H^*_B\to H^{*+4}_B$ the multiplication by the Euler class.

\begin{definition}
We say that the action $\Phi$ safisfies the property $\THL{k}$ if
$$
L^{n-k}\colon H_B^{2k}\longrightarrow H_B^{4n-2k}
$$
defined as $L^{n-k}([\alpha])=L\circ\cdots\circ L([\alpha])=[\alpha\wedge \Omega^{n-k}]$, is an isomorphism.

We say that the action $\Phi$ safisfies the {\em transversal hard Lefschetz property} (THL) if it satisfies $\THL{k}$ for every $k=0,1,\dots,n$. 
\end{definition}

\begin{definition}
The Primitive basic cohomology group of dimension $2k$ is
$$
\PHB{2k}=\set{[\beta]\in H_B^{2k} | [\beta\wedge\Omega^{n-k+1}]=0}
$$
\end{definition}

The following definition is a slight adaptation of the dimensions of \cite[Definition~1.4]{leches} to this new context.

\begin{definition}
	We will say that $\Phi$ satisfies the $k$-th {\em primitive condition} (and denote it by $P_{2k}$) if the inclusion of forms induces the following two isomorphisms:
	\begin{itemize}
		\item[$(P_1)_k$:]  $\xymatrix@1{i_{2k}=\iota_{2k}\vert_{\PHB{2k}}\colon \PHB{2k}\ar[r]^(0.7){\cong}& \HM{2k}}$;
		\item[$(P_2)_k$:]  $\HB{2k}=\PHB{2k}\oplus L(\HB{2k-4})$.
	\end{itemize}
	
\end{definition}

\begin{definition}
We say that the action $\Phi$ safisfies the property $\HL{k}$ if there exists an isomorphism
$\LL^{n-k}\colon H_M^{2k}\longrightarrow H_M^{4n-2k+3}$ such that the following diagram is commutative:
$$
\xymatrix{
          & H_M^{2k}\ar[r]^{\LL^{n-k}}   &  H_M^{4n-2k+3}\ar[d]^{\rho} \\
\PHB{2k}\ar[ru]^{i_{2k}}\ar@{^{(}->}[r] & H_B^{2k}\ar[r]^{L^{n-k}}   &  H_B^{4n-2k}},
$$
where $i_k$ is induced by the inclusion of forms.

We say that the action $\Phi$ safisfies the {\em hard Lefschetz property} (THL) if it satisfies $\HL{k}$ for every $k=0,1,\dots,n$.
\end{definition}

For a certain kind of $\ese^3$-actions both $\HL{}$ and $\THL{}$ will be equivalent.

\begin{definition}
We say that the almost free action $\Phi$ satisfies the {\em injectivity condition} if the multiplication by the Euler class induces a monomorphism:
$
L\colon \HB{2k-3} \rightarrowtail \HB{2k+1}
$
for $k=0,1,\dots,n$.
\end{definition}
Notice that this condition is not present in the context of flows as a necessary hypothesis for the equivalence of $\THL{}$ and $\HL{}$. The reason is that for $ \ese^3$-actions $\THL{k-1}$ is not enough to guarantee the necessary Primitive condition $P_{2k}$ which is the key to get $\HL{k}$. We have

\begin{proposition}\label{prop:pk} For every $k\le n$,
	$$
	\left.
	\begin{matrix}
		L\colon H_B^{2k-3}\longrightarrow H_B ^{2k+1} \text{monomorphism}\\
		\text{and}\\
		(THL)_{k-2}
	\end{matrix}
	\right\}
	\Longrightarrow P_{2k}
	$$
\end{proposition}

\begin{proof}
The monomorphism of the statement and the Gysin sequence \eqref{eq:gysin} yield immediately that 
$
\iota_{2k}\colon\HB{2k}\longrightarrow\HM{2k}
$
is an epimorphism. Now, the proof of the  result \cite[Proposition~1.8]{leches}, which is analogous to this, applies straightforwardly, having on account that the epimorphism we just got is the same property one gets in \cite{leches} from $\THL{k-1}$ and \cite[Lemma~1.6]{leches} at the begining of that proof, and of course, changing $e$ by $\Omega$.
\end{proof}

Our main result is the following:

\begin{theorem}\label{th:equiv}
Let $\Phi$ be an almost free smooth action of $\ese^3$ on the closed manifold $M$ satisfying the injectivity condition. Then, it satisfies $\HL{}$ if and only if it satisfies $\THL{}$.
\end{theorem}

\begin{proof}
In fact, we will prove something slightly stronger: given $k\le n$, provided we have monomorphisms 
$
L\colon \HB{2j-3} \rightarrowtail \HB{2j+1}
$ for all $j\le k$, then we have $\THL{j} \Longleftrightarrow \HL{j}$ for all $j\le k$. 
Again, the proof of the analogous \cite[Theorem~1.10]{leches} applies straightforwardly, having on account for the reciprocal implication that the injectivity condition ensures the first hypothesis of Proposition~\ref{prop:pk}, eventually yielding the needed isomorphisms of $P_{2k}$.
\end{proof}

This last result justifies the following definition.

\begin{definition}
	We say that an almost-free action $\Phi$ of $\ese^3$ on a closed manifold $M$ is a {\em Lefschetz $\ese^3$-action} if it satisfies the injectivity condition and $\THL{}$ or $\HL{}$.
\end{definition}

\section{3-Sasakian manifolds and other examples}

Our next goal is to show that 3-Sasakian manifolds satisfy both $\HL{}$ and $\THL{}$. Recall (see \cite[chapter~13]{galicki} for the basic facts about 3-Sasakian structures) that a 3-Sasakian structure (which has associated three fundamental vector fields and an invariant metric making them orthonormal, so that each of the vector fields is the Reeb vector field of a Sasakian structure) over a manifold of dimension $4n+3$ has associated an almost-free action of $\ese^3$ for which the metric is $\ese^3$-invariant and the vector fields are the fundamental vector fields of the action (up to convention, see Remark \ref{rem:su2}). When that action is, actually, free,  the 3-Sasakian structure is said to be {\em regular}.

It is known that if $\Phi$ is the action of the 3-sphere induced by a 3-Sasakian structure, then we all the small basic Betti numbers are zero, that is, 
\begin{equation}\label{eq:smallbetti}
b_B^{2k+1}=0,\ \forall k=0,1,\dots,n
\end{equation}
and so, all 3-Sasakian manifolds satisfy the injectivity condition, trivially. By Theorem \ref{th:equiv} if we show that 3-Sasakians satisfy one of the defined Hard Lefschetz properties, they will satisfy both of them.

To prove that 3-Sasakian manifolds are Lefschetz we shall need another Gysin-like sequence associated to any almost-free action $\Phi$. Consider one of the fundamental vector fields of $\Phi$, say $X_1$, and its corresponding isometric flow, which is an almost-free $\ese^1$-action. We consider the complex  $\Omega_Z^*$ of basic forms of $M$ with respect to that foliation by circles and its cohomology $H_Z^*$. Notice that any form basic for the action of $\ese^3$ is automatically basic for the action of $\ese^1$. Thus, we have the following short exact sequence of complexes:
$$
0\to \Omega_B^*\to \Omega_Z^* \to K^*(M)=\frac{\Omega^*_Z}{\Omega^*_B}\to 0
$$
This sequence induces a long exact sequence in cohomology of Gysin type. The following result determines the cohomology of the so-called {\em Gysin term} $K^*(M)$:
\begin{proposition}\label{prop:s2}
In the previous conditions, the map 
$$
f\colon\Omega^*_B\longrightarrow K^{*+2}(M),
$$
defined by $f(\omega)=\overline{d\chi_1\wedge\omega}$ induces an isomorphism in cohomology.
\end{proposition}

\begin{proof}
It is easily checked that the map is well defined, commutes with the differential and, moreover, it does not depend on the invariant metric chosen. By Mayer-Vietoris, we can restrict ourselves to prove the statement by taking $M$ to be the tubular neighbourhood of a stratum $S$ of the action $\Phi$ (for the basics of compact transformation groups we refer the reader to \cite{bredon72}), and by equivariant retraction, to a stratum $S$. We have that any stratum is a twisted product by the normalizer of a certain finite isotropy group $H$:
$$
S=\ese^3\times_{N(H)}S^H,
$$
where $S^H$ is the submanifold of $S$ formed by the points whose isotropy subgroup is $H$. As $N(H)/H$ acts freely on $S^H$ we have 
$$
S=\ese^3\times_{N(H)}(\R^m\times N(H)/H),
$$
so by equivariant retraction we can restrict ourselves to the case
$$
S=\ese^3\times_{N(H)}(N(H)/H)\cong \ese^3/H
$$
It remains to address the case $M=\ese^3/H$. As $B=M/\ese^3$ is a point, we have
$$
H^0(K(M))=0,\quad\text{and}\quad H^k(K(M))=H^k_Z,\quad\forall k>0.
$$
By putting $\ese^2=\ese^3/\ese^1$, we have that 
$$
\Omega_Z^*=\Omega^*(\ese^2/H)=\Omega^*(\ese^2)^H.
$$
As $H$ is a finite group, we finally obtain
$$
H^*_Z=H(\Omega^*(\ese^2)^H)=H^*(\ese^2)^H=H^*(\ese^2)=H^2(\ese^2)=\R
$$
So, $f$ induces an isomorphism in cohomology and we are done.
\end{proof}

From Proposition \ref{prop:s2} we get the following Gysin seguence, which we can informally refer to as the Gysin sequence {\em of the $\ese^2$-action}.
\begin{corollary}\label{cor:s2}
Let $\Phi$ be an almost-free action of $\ese^3$ on a closed manifold, consider $\ese^1$ as a subgroup of $\ese^3$ and $\Phi_1$ its corresponding action on $M$. Then, we have a long exact sequence
$$
\dots\to  H^k_B\longrightarrow H^k_Z\longrightarrow H^{k-2}_B\longrightarrow H^{k+1}_B\to\dots 
$$
where $H^*_B$ and $H^*_Z$ stand for the basic cohomology of $\Phi$ and $\Phi_1$, respectively.
\end{corollary}
Notice that there is no orbifold cohomology involved in the previous sequence: all the cohomology groups are defined by using differential forms defined in $M$. Now, we have all the ingredients to show that all 3-Sasakian manifolds are Lefschetz.

\begin{theorem}
	Let $\Phi$ be the action of $\ese^3$ associated to a 3-Sasakian structure over a closed manifold $M$. Then, $\Phi$ satisfies both Hard Lefschetz Properties.
\end{theorem}

\begin{proof}
	 Consider one of the Sasakian structures (say $X_1$) and its corresponding isometric flow, which is indeed an almost-free $\ese^1$-action. We have
	 \begin{equation}\label{eq:chichi}
	 \begin{split}
	 d\chi_1\wedge d\chi_1 &= d(\chi_1\wedge d\chi_1) = d(\chi_1\wedge (e_1-\chi_2\wedge\chi_3))\\
	 	&=(e_1^2 + e_2^2 + e_3^2) - d(e_2\wedge\chi_2 + e_3\wedge\chi_3),
	 \end{split}
	 \end{equation}
	 where we have used \eqref{eq:omega}.
	 Now, notice that the forms $d\chi_1\wedge d\chi_1, \Omega$ and $e_2\wedge\chi_2 + e_3\wedge\chi_3$ are $X_1$-basic. So, from \eqref{eq:chichi} we get
	$$
	[d\chi_1\wedge d\chi_1]=[\Omega]\in H^4(M/\F_1)
	$$
	and, correspondingly, we have a commutative diagram
	$$
	\xymatrix@C=40pt{H_Z^{2k} \ar[r]^{L_e^{2n-2k}}& H_Z^{4n-2k}\\
		H_B^{2k} \ar[r]^{L_{\Omega}^{n-k}}
		\ar@{>->}[u] & H_B^{4n-2k}\ar@{>->}[u]}
	$$
	where $H_Z^i=H^i(M/\F_1)$ and $L_e$ and $L_{\Omega}$ stand for the multiplication by the corresponding Euler classes $[e_1]$ and $[\Omega]$, respectively. 
	The injectivity of the vertical arrows follows from Corollary \ref{cor:s2} and \eqref{eq:smallbetti}. Then, as the bottom groups of the diagram are of the same dimension and the top horizontal arrow is an isomorphism (because $\F_1$ is a Sasakian manifold, and thus, it is a Lefschetz flow and satisfies $\THL{}$ for flows), we get that the bottom horizontal arrow is an isomorphism, and thus we get $\THL{}$ for $\Phi$ and, by virtue of Theorem \ref{th:equiv}, $\HL{}$.
\end{proof}

\begin{example} We finish this paper by showing that, apart from the ones arising from 3-Sasakian structures, there are many more $\ese^3$-actions which satisfy both versions of the HLP. On one hand, in \cite{hsiang66} it is shown that there exist infinitely many differentiably inequivalent free actions of $\ese^3$ on homotopy spheres of dimension $4n+3$ for every integer $n\ge2$. Should any of these actions come from a 3-Sasakian manifold, then it must be regular, because they are free.  On the other hand, by \cite[Proposition~13.5.6]{galicki}, given $n\ge1$, there are only finitely many regular 3-Sasakian manifolds of dimension $4n+3$, up to isometries. As a consequence, given $4n+3$ with $n\ge 2$, we can chose an action (in fact, infinitely many unequivalent ones) on a homotopy sphere $M$ that doesn't come from a 3-Sasakian structure. As $M$ has the cohomology of $\ese^{4n+3}$, the Gysin sequence \eqref{eq:gysin}  immediately yields, along with $H_B^0=\R$, that the connecting morphisms
$$
H_B^{4i}\stackrel{\wedge[\Omega]}{\longrightarrow} H_B^{4(i+1)},\quad i=0,1,\dots,n-1
$$
are isomorphisms and that the remaining basic cohomology groups  are trivial. So, the basic cohomology is isomorphic to the cohomology of the quaternionic projective space $\mathbb{H}P^n$, and satisfies both the injectivity condition and the transversal Lefschetz Property, hence, by Theorem \ref{th:equiv}, also the $\HL{}$.

Finally, notice that we cannot, in general, assure whether we can pick $M$ to be the standard sphere $\ese^{4n+3}$, but for dimension 11, we actually can, because in \cite{hsiangs11} the author constructs infinitely many differentiably unequivalent free actions of $\ese^3$ on $\ese^{11}$.

\end{example}
	
%\bibliography{nirebib}

%\bibliographystyle{alpha}
%\bibliographystyle{unsrt}
%\bibliographystyle{abbrv}
%\bibliographystyle{plain}
%\bibliographystyle{smfalpha}
%\bibliographystyle{smfplain}
%\bibliographystyle{amsplain}
%\bibliographystyle{apa}

\printbibliography

\end{document}